\documentclass[dvips,coll,authoryear]{imsart}
\pdfoutput=1

\usepackage{amsmath, amsthm}
\usepackage[pdftex]{graphicx}

\usepackage{amssymb}
\usepackage{pifont}

\volumetitle{???}
\volume{0}
\pubyear{0000}
\firstpage{1}
\lastpage{1}

\startlocaldefs

\def\ni{\noindent}

\DeclareMathOperator*{\diag}{diag}
\DeclareMathOperator{\Var}{Var}
\DeclareMathOperator{\tr}{tr}

\def\pa{\partial}
\def\para{\parallel}
\def\tsum{\textstyle{\sum}}

\def\al{\alpha}
\def\be{\beta}
\def\hbe{\hat{\beta}}
\def\de{\delta}
\def\ep{\epsilon}
\def\ga{\gamma}
\def\hga{\hat{\gamma}}

\def\la{\lambda}

\def\hmu{\hat{\mu}}
\def\si{\sigma}

\def\ze{\zeta}

\def\bR{\mathbb{R}}
\def\bI{\mathbb{I}}

\def\cA{{\cal A}}
\def\chA{\hat{\cal A}}
\def\cB{{\cal B}}
\def\cC{{\cal C}}
\def\cD{{\cal D}}
\def\cE{{\cal E}}
\def\hg{\hat{g}}
\def\cH{{\cal H}}
\def\cO{{\cal O}}

\def\hq{\hat{q}}
\def\tX{\tilde{X}}
\def\bX{\bar{X}}
\def\dX{\ddot{X}}
\def\tY{\tilde{Y}}
\def\bY{\bar{Y}}

\def\etal{{et al.}~}
\def\etalc{et al.,~}

\def\beqn{\begin{equation}}
\def\neqn{\end{equation}}

\newcommand\T{\rule{0pt}{2.6ex}}
\newcommand\B{\rule[-1.2ex]{0pt}{0pt}}

\newcommand{\tick}{\ding{52}}
\newcommand{\cross}{\ding{56}}

\numberwithin{equation}{section}
\theoremstyle{plain}
\newtheorem{thm}{Theorem}[section]
\newtheorem{lmm}[thm]{Lemma}

\newtheorem{defn}{Definition}

\endlocaldefs

\begin{document}

\begin{frontmatter}

\title{Local polynomial regression and variable selection}

\runtitle{Local variable selection}

\begin{aug}
\author{\fnms{Hugh} \snm{Miller}\ead[label=e1]{h.miller@ms.unimelb.edu.au}}
\and
\author{\fnms{Peter} \snm{Hall}\ead[label=e2]{halpstat@ms.unimelb.edu.au}}

\runauthor{H. Miller and P. Hall}

\affiliation{The University of Melbourne}

\address[a]{Department of Mathematics and Statistics, The University of Melbourne, Melbourne, VIC 3052, Australia, \printead{e1,e2}}

\end{aug}

\begin{abstract}
We propose a method for incorporating variable selection into local polynomial regression. This can improve the accuracy of the regression by extending the bandwidth in directions corresponding to those variables judged to be are unimportant. It also increases our understanding of the dataset by highlighting areas where these variables are redundant. The approach has the potential to effect complete variable removal as well as perform partial removal when a variable redundancy applies only to particular regions of the data. We define a nonparametric oracle property and show that this is more than satisfied by our approach under asymptotic analysis. The usefulness of the method is demonstrated through simulated and real data numerical examples.
\end{abstract}



\end{frontmatter}

\section{Introduction}

The classical regression problem is concerned with predicting a noisy continuous response using a $d$-dimensional predictor vector with support on some $d$-dimensional subspace. This functional relationship is often taken to be smooth and methods for estimating it range from parametric models, which specify the form of the relationship between predictors and response, through to nonparametric models, which have fewer prior assumptions about the shape of the fit. An important consideration for fitting such a regression model is whether all $d$ predictors are in fact necessary. If a particular predictor has no relationship to the response, the model be made both simpler and more accurate by removing it. In recent years there has been strong interest in techniques that automatically generate such ``sparse'' models. Most attention has been given to parametric forms, and in particular the linear model, where the response is assumed to vary linearly with the predictors. There has also been some investigation into variable selection for nonlinear models, notably through the use of smoothing splines and local regression.

One common feature of the existing sparse methods is that the variable selection is ``global'' in nature, attempting to universally include or exclude a predictor. Such an approach does not naturally reconcile well with some nonparametric techniques, such as local polynomial regression, which focus on a ``local'' subset of the data to estimate the response. In this local context it would be more helpful to understand local variable influence, since predictors that are irrelevant in some regions may in fact be important elsewhere in the subspace. Just as in the global setting, such information would allow us to improve the accuracy and parsimony of a model, but at a local level.

However this approach to variable selection can be problematic. Most notably, variable significance  affects the definition of ``local''. To illustrate concretely, suppose that two data points were close in every dimension except one. In typical local regression these points would not be considered close, and so the response at one point would not impact the other. If, however, we establish that the one predictor they differ by is not influential over a range that includes both these points, then they should actually be treated as neighbouring, and be treated as such in the model. Any methodology seeking to incorporate local variable influence needs to accommodate such potential situations.

Understanding local variable significance can also give additional insight into a dataset. If a variable is not important in certain regions of the support, knowledge of this allows us to discount it in certain circumstances, simplifying our understanding of the problem. For example, if none of the variables are relevant in a region, we may treat the response as locally constant and so know that we can ignore predictor effects when an observation lies in this region.

A final consideration is theoretical performance. In particular we shall present on approach that is ``oracle''; that is, its performance is comparable to that of a particularly well-informed statistician, who has been provided in advance with the correct variables. It is interesting to note that variable interactions often cause sparse parametric approaches to fail to be oracle, but in the local nonparametric setting this is not an issue, because such interactions vanish as the neighbourhood of consideration shrinks.

In this paper we propose a flexible and adaptive approach to local variable selection using local polynomial regression. The key technique is careful adjustment of the local regression bandwidths to allow for variable redundancy. The method has been named LABAVS, standing for ``locally adaptive bandwidth and variable selection''. Section~2 will introduce the LABAVS algorithm, including a motivating example and possible variations. Section~3 will deal with theoretical properties and in particular establishes a result showing that the performance of LABAVS is better than oracle when the dimension remains fixed. Section~4 presents numerical results for both real and simulated data, showing that the algorithm can improve prediction accuracy and is also a useful tool in arriving at an intuitive understanding of the data. Technical details have been relegated to an appendix which may be found in the long version of this paper (Miller and Hall, 2010).

LABAVS is perhaps best viewed as an improvement to local polynomial regression, and will retain some of the advantages and disadvantages associated with this approach. In particular, it still suffers the ``curse of dimensionality,'' in that it struggles to detect local patterns when the dimension of genuine variables increases beyond a few. It is not the first attempt at incorporating variable selection into local polynomial regression; the papers by Lafferty and Wasserman (2008) and Bertin and Lecu\'{e} (2008) also do this. We compare our approach to these in some detail in Section~\ref{se:compare}. LABAVS can also be compared to other nonparametric techniques in use for low to moderate dimensions. These include generalised additive models, MARS and tree based methods (see Hastie \etalc 2001).

The earliest work on local polynomial regression dates back to that of Nadaraya (1964) and Watson (1964). General references on the subject include Wand and Jones (1995), Simonoff (1996) and Loader (1999). An adaptive approach to bandwidth selection may be found in Fan and Gijbels (1995), although this was not in the context of variable selection. Tibshirani (1996) studies the LASSO, one of the most popular sparse solutions for the linear model; more recent related work on the linear model includes that of Candes and Tao (2007) and Bickel \etal (2009). Zou (2006) created the adaptive version of the LASSO and proved oracle performance for it. Lin and Zhang (2006) and Yuan and Lin (2006) have investigated sparse solutions to smoothing spline models.

The LABAVS algorithm also bears some similarity to the approach adopted by Hall \etal (2004). There the aim was to estimate the conditional density of a response using the predictors. Cross-validation was employed and the bandwidths in irrelevant dimensions diverged, thereby greatly downweighting those components. In the present paper the focus is more explicitly on variable selection, as well as attempting to capture local variable dependencies.

\section{Model and Methodology}

\subsection{Model and definitions}

Suppose that we have a continuous response $Y_i$ and a $d$-dimensional random predictor vector $X_i = (X^{(1)}_i, \ldots, X^{(d)}_i)$ which has support on some subspace $\cC \subset \bR^d$. Further, assume that the observation pairs $(Y_i,X_i)$ are independent and identically distributed for $i=1,\ldots,n$, and that $X_i$ has density function $f$. The response is related to the predictors through a function $g$,
\beqn Y_i = g(X_i) + \ep_i \,, \label{eq:simple} \neqn
with the error $\ep_i$ having zero mean and fixed variance. Smoothness conditions for $f$ and $g$ will be discussed in the theory section.

Local polynomial regression makes use of a kernel and bandwidth to assign increased weight to neighbouring observations compared to those further away, which will often have zero weight. We take $K(u) = \prod_{1 \leq j \leq d} K^*(u^{(j)})$ to be the $d$-dimensional rectangular kernel formed from a one dimensional kernel $K^*$ such as the tricubic kernel,
$$ K^*(u^{(j)})=(35/32)(1-x^2)^3 \bI(|x|<1) \,.$$
\label{kernel} Assume $K^*$ is symmetric with support on $[-1,1]$. For $d \times d$ bandwidth matrix $H$ the kernel with bandwidth $H$, denoted $K_H$, is
\beqn K_H(u) = \frac{1}{|H|^{1/2}}K(H^{-1/2}u) \,. \label{eq:ker} \neqn
We assume that the bandwidth matrices are diagonal, $H = \diag(h_1^2,\ldots,h_d^2)$, with each $h_j>0$, and write $H(x)$ when $H$ varies as a function of $x$. Asymmetric bandwidths can be defined as having both a lower and an upper (diagonal) bandwidth matrix,  $H^L$ and $H^U$ respectively, for a given estimation point $x$, rather than a single bandwidth $H$ for all $x$. The kernel weight of an observation $X_i$ at estimation point $x$ with asymmetrical bandwidth matrices $H^L(x)$ and $H^U(x)$, is
\begin{eqnarray}
K_{H^U(x),H^L(x)}(X_i - x) &=& \prod_{j \,:\, X_i^{(j)} < x^{(j)}} \frac{1}{h^L_j(x)} K^*\left( \frac{X^{(j)}_i - x^{(j)}}{h^L_j(x)} \right) \\
&& \qquad \cdot \prod_{j\,:\, X_i^{(j)} \ge x^{(j)}} \frac{1}{h^U_j(x)} K^*\left( \frac{X^{(j)}_i - x^{(j)}}{h^U_j(x)} \right) \,. \nonumber
\end{eqnarray}
This amounts to having (possibly) different window sizes above and below $x$ in each direction. Although such unbalanced bandwidths would often lead to undesirable bias properties in local regression, here they will be used principally to extend bandwidths in dimensions considered redundant, so this issue is not a concern.

We also allow the possibility of infinite bandwidths $h_j = \infty$. In calculating the kernel in (\ref{eq:ker}) when $h_j$ is infinite, proceed as if the $j$th dimension did not exist (or equivalently, as if the $j$th factor in rectangular kernel product is always equal to 1). If all bandwidths are infinite, consider the kernel weight to be 1 everywhere. Although the kernel and bandwidth conditions above have been defined fairly narrowly to promote simplicity in exposition, many of these assumptions are easily generalised.

Local polynomial regression estimates of the response at point $x$, $\hg(x)$, are found by fitting a polynomial $q$ to the observed data, using the kernel and bandwidth to weight observations. This is usually done by minimising the weighted sum of squares,
\beqn \sum_{i=1}^n \{Y_i - q(X_i-x)\}^2 K_H(X_i-x) \,. \label{eq:locpoly} \neqn
Once the minimisation has been performed, $q(0)$ becomes the point estimate for $g(x)$. The polynomial is of some fixed degree $p$, with larger values of $p$ generally decreasing bias at the cost of increased variance. Of particular interest in the theoretical section will be the local linear fit, which minimises
\beqn  \sum_{i=1}^n \left[ \left\{ Y_i - \ga_0 - \sum_{j=1}^n (X_i^{(j)}-x^{(j)})\ga_j   \right\}^2 K_{H}(X_i-x) \right] \,,\label{eq:loclin} \neqn
over $\ga_0$ and $\ga=(\ga_1,\ldots,\ga_d)$.

\subsection{The LABAVS Algorithm}

Below is the LABAVS algorithm that will perform local variable selection and vary the bandwidths accordingly. The choice of $H$ in the first step can be local or global and should be selected as for a traditional polynomial regression, using cross-validation, a plug-in estimator or some other standard technique. Methods for assessing variable significance in Step~2, and the degree of shrinkage needed in Step~4, are discussed below.

$$\begin{array}{p{13cm}}
\ni {\sl LABAVS Algorithm}
\begin{enumerate}
\item Find a starting $d\times d$ bandwidth $H = \diag(h^2,\ldots,h^2)$.
\item For each point $x$ of a representative grid in the data support, perform local variable selection to determine disjoint index sets $\chA^+(x),\chA^-(x)$, with $\chA^+(x) \cup \chA^-(x) = \{1,\ldots,d\}$,  for variables that are considered relevant and redundant respectively.
\item For any given $x$, derive new local bandwidth matrices $H^L(x)$ and $H^U(x)$ by extending the bandwidth in each dimension indexed in $\chA^-(x)$. The resulting space given nonzero weight by the kernel $K_{H^L(x),H^U(x)}(u-x)$ is the rectangle of maximal area with all grid points $x_0$ inside the rectangle satisfying $\chA^+(x_0) \subset \chA^+(x)$. Here $\chA^+(x)$ is calculated explicitly as in Step~2, or taken as the set corresponding the closest grid point to $x$.
\item Shrink the bandwidth slightly for those variables in $\chA^+(x)$ according to the amount that bandwidths have increased in the other variables.
\item Compute the local polynomial estimator at $x$, excluding variables in $\chA^-(x)$ and using adjusted assymetrical bandwidths $H^L(x)$ and $H^U(x)$. The expression to be minimised is
$$ \sum_{i=1}^n \left\{ Y_i - q(X_i-x)   \right\}^2 K_{H^L(x), H^U(x)}(X_i-x)  \,, $$
where the minimisation runs over all polynomials $q$ of appropriate degree. The value of $q(0)$ in the minimisation is the final local linear estimator.
\end{enumerate}
\end{array}$$

The key feature of the algorithm is that variable selection directly affects the bandwidth, increasing it in the direction of variables that have no influence on the point estimator. If a variable has no influence anywhere, it has the potential to be completely removed from the local regression, reducing the dimension of the problem. For variables that have no influence in certain areas, the algorithm achieves a partial dimension reduction. The increased bandwidths reduce the variance of the estimate and Step~4 swaps some of this reduction for a decrease in the bias to further improve the overall estimator.

\label{se:eg}

As a concrete example of the approach, define the following one-dimensional ``huberised'' linear function:
\beqn g(x) = x^2\bI(0 < x \le 0.4) + (0.8x-0.16)\bI(x>0.4) \,,\neqn
and let $g(X) = g\{([X^{(1)}]_+^2+[X^{(2)}]_+^2)^{1/2}\}$ for 2-dimensional random variable $X = (X^{(1)}, X^{(2)})$. Assume that $X$ is uniformly distributed on the space $[-2,2]\times[-2,2]$. Notice that when $X^{(1)}$, $X^{(2)} <0$ the response variable $Y$ in (\ref{eq:simple}) is independent of $X^{(1)}$ and $X^{(2)}$; when $X^{(1)} < 0$ and $X^{(2)} >0$ the response depends on  $X^{(2)}$ only; when $X^{(1)} > 0$ and $X^{(2)} <0$ the response depends on  $X^{(1)}$ only; when $X^{(1)}$, $X^{(2)} >0$ the response depends on both $X^{(1)}$ and $X^{(2)}$. Thus in each of these quadrants a different subset of the predictors is significant. A local approach to variable significance can capture these different dependencies, while a global variable redundancy test would not eliminate any variables.

\begin{figure}[htbp] \begin{center}
\includegraphics[width = 3.5in, height=3.5in]{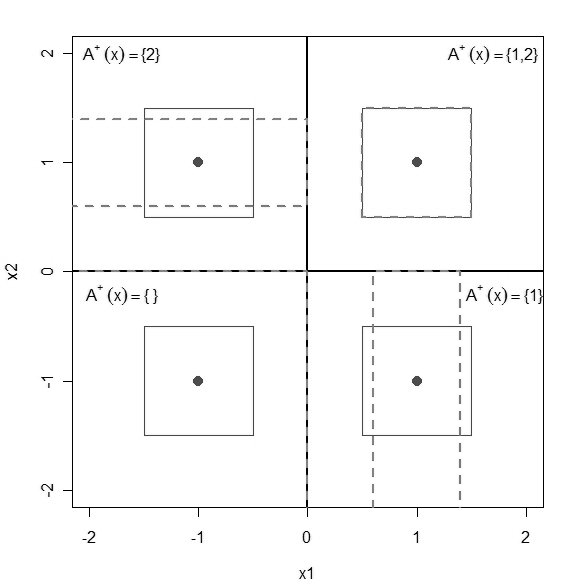}
\vspace{-8pt}
{\it \caption{ Bandwidth adjustments under ideal circumstances in illustrative example. \label{fig2}}}
\end{center} \end{figure}

Now consider how the algorithm applies to this example, starting with a uniform initial bandwidth of $h=0.5$ in both dimensions. Assuming that variable significance perfectly on a dense grid, figure~\ref{fig2} illustrates the adjusted bandwidths for each of the quadrants. The dots are four sample estimation points, the surrounding unit squares indicate the initial bandwidths and the dashed lines indicate how the bandwidths are modified. In the bottom left quadrant both variables are considered redundant, and so the bandwidth expands to cover the entire quadrant. This is optimal behaviour, since the true function is constant over this region, implying that the best estimator will be produced by including the whole area. In the bottom right quadrant the first dimension is significant while the second is not. Thus the bandwidth for the second dimension is ``stretched'', while the first is shrunken somewhat. Again, this is desirable for improving the estimator. The stretching in the second dimension improves the estimator by reducing the variance as more points are considered. Then the shrunken first dimension swaps some of this reduction in variance for decreased bias. Finally, in the top right quadrant, there is no change in the bandwidth since both variables are considered to be significant.

\subsection{Variable selection step} \label{se:approaches}
Below are three possible ways to effect variable selection at $x_0$ in Step~2 of the algorithm, presented in the context of local linear regression. They all make use of a tuning parameter $\la$ which controls how aggressive the model is in declaring variables as irrelevant. Cross validation can be used to select an appropriate level for $\la$. So that the tuning parameters are comparable at different points in the data domain, it is useful to consider a local standardisation of the data at $x_0$. Define $\bX_{x_0} = (\bX_{x_0}^{(1)}, \ldots, \bX_{x_0}^{(d)})$ and $\bY_{x_0}$ by
\begin{eqnarray*}
 \bX_{x_0}^{(j)} &=& \frac{\sum_{i=1}^n X_i^{(j)} K_H(X_i-x_0)}{ \sum_{i=1}^n K_H(X_i-x_0)} \,, \quad \bY_{x_0} = \frac{\sum_{i=1}^n Y_i K_H(X_i-x_0)}{ \sum_{i=1}^n K_H(X_i-x_0)} \,,
\end{eqnarray*}
 and define $\tX_i = (\tX_i^{(1)},\ldots,\tX_i^{(d)})$ and $\tY_i$ by
\begin{eqnarray}
 \tX_i^{(j)} &=& \frac{(X^{(j)}_i-\bX^{(j)}_{x_0})K_H(X_i-x_0)^{1/2}} {\big\{\sum_{i=1}^n \big(X_i^{(j)} - \bX^{(j)}_{x_0} \big)^2 K_H(X_i-x_0) \big\}^{1/2}} \,, \label{eq:scale} \\
 \tY_i &=& (Y_i-\bY_{x_0})K_H(X_i-x_0)^{1/2} \,. \nonumber
\end{eqnarray}
Notice that $\tX$ and $\tY$ incorporate the weight $K_H(X_i-x_0)$ into the expression.

\begin{enumerate}
\item {\bf Hard thresholding:} Choose parameters to minimise the weighted least squares expression,
\beqn \sum_{i=1}^n \left\{ \tY_i - \be_0 - \sum_{j =1}^n \tX_i^{(j)}\be_j   \right\}^2  \,, \label{eq:appr1} \neqn
and classify as redundant those variables for which $|\hbe_j|<\la$. This can be extended to higher degree polynomials, although performance tends to be more unstable.
\item {\bf Backwards stepwise approach:} For each individual $j$, calculate the percentage increase in the sum of squares if the $j$th variable is excluded from the local fit. Explicitly, if $\hq$ is the optimal local fit using all variables and $\hq_j$ is the fit using all except the $j$th, we classify the $j$th variable as redundant if
\beqn \frac{\sum_{i=1}^n \left\{ Y_i - \hq_j(\hat{X}_i)   \right\}^2 K_H(\hat{X}_i) - \sum_{i=1}^n \left\{ Y_i - \hq(\hat{X}_i)   \right\}^2 K_H(\hat{X}_i)}{\sum_{i=1}^n \left\{ Y_i - \hq(\hat{X}_i)   \right\}^2 K_H(\hat{X}_i)} < \la \,, \label{eq:appr2} \neqn
where $\hat{X}_i =X_i-x_0$.
\item {\bf Local lasso:} Minimise the expression
\beqn \sum_{i=1}^n  \left\{ \tY_i - \ga_0 - \sum_{j =1}^n \tX_i^{(j)}\ga_j   \right\}^2  + \la \sum_{j=1}^d|\ga_j| \,. \label{eq:appr3} \neqn
Those variables for which $\ga_j$ are set to zero in this minimisation are then classed as redundant. While the normal lasso can have consistency problems (Zou, 2006), this local version does not since variables are asymptotically independent as $h\to 0$. The approach also scales naturally to higher order polynomials, provided all polynomial terms are locally standardised; a variable is considered redundant if all terms that include it have corresponding parameters set to zero by the lasso.
\end{enumerate}

We have found that the first and second of the above approaches have produced the most compelling numerical results. The numerical work in Section~\ref{se:num} uses the first approach for linear polynomials, while the theoretical work in Section~\ref{theory} establishes uniform consistency for both of the first two methods, guaranteeing oracle performance.

\subsection{Variable shrinkage step}

The variable shrinkage step depends on whether the initial bandwidth, and thus the shrunken bandwidth $h'$, is chosen locally or globally. Define
\beqn V(x,H) = \frac{\tsum K_{H}(X_i-x)^2 }{ \{\tsum K_{H}(X_i-x) \}^2 } \,, \neqn
where the bandwidth term in the function $V$ is allowed to be asymmetrical. Then in the local case, letting $d'(x)$ denote the cardinality of $\chA(x)$, let
\beqn M(x) = V[x,\{H^L(x),H^L(x)\}] / V(x,H)  \label{eq:fac} \,. \neqn
The expression is asymptotically proportional to $\{h'(x)\}^{-d'(x)}$ and estimates the degree of variance stabilisation resulting from the bandwidth adjustment. Using this, the correct amount of bandwidth needed in step~4 is $h'(x) = h \{M(x) d'(x)/d\}^{1/4} $. Since both sides of this expression depend on $h'(x)$, shrinkage can be approximated in the following way. Let
$$ M^*(x) = V[x,\{\tilde{H}^L(x),\tilde{H}^L(x)\}] / V(x,H) \,,$$
where $\tilde{H}^L(x)$ and $\tilde{H}^U(x)$ are the bandwidth matrices immediately after step 3. Then the shrunken bandwidths are $h'(x) = h\{ M^*(x) d'(x)/d )^{1/(d'(x)+4)}$.

In the global bandwidth case, we define
\beqn M[\{H^L(X),H^U(X)\},H] = \frac{E\left(V[X,\{H^L(X),H^L(X)\}]\right)}{E\{V(X,H)\}} \,. \label{newfac1} \neqn
This expression measures the average variance stabilisation across the domain. In this case, the shrinkage factor should satisfy
\beqn h' = h \left( M[\{H^L(X),H^U(X)\},H] E\{d'(X)\}/d \right)^{1/4} \,. \label{newfac2} \neqn
The theoretical properties in Section~\ref{theory} deal with the global bandwidth scenario. The treatment for the local case is similar, except that care must be taken in regions of the domain where the function $g$ behaves in a way that is exactly estimable by a local polynomial and thus has potentially no bias.

\subsection{Further remarks}

\begin{enumerate}
\item The choice of distance between grid points in Step~2 is somewhat arbitrary, but should be taken as less than $h$ so that all data points are considered in calculations. In the asymptotic theory we let this length decrease faster that the rate of the bandwidth, and in numerical experimentation the choice impacts only slightly on the results.
\item Step 5 of the algorithm forces the estimate at point $x$ to exclude variables indexed in $\chA^-(x)$. An alternative is to still use all variables in the final fit. This may be advantageous in situations with significant noise, where variable admission and omission is more likely to have errors. Despite including these extra variables, the adjusted bandwidths still ensure that estimation accuracy is increased.
\item Finding the maximal rectangle for each representative point, as suggested in step~3 of the algorithm, can be a fairly intensive computational task. In our numerical work we simplified this by expanding the rectangle equally until the boundary met a ``bad'' grid point (i.e.~a point $x'$ such that $\chA^+(x')\nsubseteq \chA^+(x)$). The corresponding direction was then held constant while the others continue to increase uniformly. We continued until each dimension stopped expanding or grew to be infinite. This approach does not invalidate the asymptotic results in Section~\ref{theory}, but there may be some deterioration in numerical performance associated with this simplification.
\item If a variable is redundant everywhere, results in Section~\ref{theory} demonstrate that the algorithm is consistent; the probability that the variable is classified as redundant everywhere tends to 1 as $n$ grows. However, the exact probability is not easy to calculate and for fixed $n$ we may want greater control over the ability to exclude a variable completely. In such circumstances a global variable selection approach may be appropriate.
\item As noted at the start of Section~2.2, the initial bandwidth in Step~1 does not necessarily have to be fixed over the domain.  For instance, a nearest neighbour bandwidth, where $h$ at $x$ is roughly proportional to $f(x)^{-1}$, could be used. Employing this approach offers many practical advantages and the theoretical basis is similar to that for the constant bandwidth. The numerical work makes use of nearest neighbour bandwidths throughout. In addition, we could use an initial bandwidth that was allowed to vary for each variable, $H = \diag(h_1^2,\ldots,h_p^2)$. So long as, asymptotically, each $h_j$ was equal to $C_j h$ for some controlling bandwidth $h$ and constant $C_j$, the theory would hold, although details are not pursued here.
\end{enumerate}

\subsection{Comparison to other local variable selection approaches} \label{se:compare}

As mentioned in the introduction, two recent papers take a similar approach to this problem. Firstly Lafferty and Wasserman (2008) introduce the rodeo procedure. This attempts to assign adaptive bandwidths based on the derivative with respect to the bandwidth for each dimension, $\partial \hg(x)/\partial h_j$. This has the attractive feature of bypassing the actual local shape and instead focussing on whether an estimate is improved by shrinking the bandwidths. It is also a greedy approach, starting with large bandwidths in each direction and shrinking only those that cause a change in the estimator at a point. The second paper is by Bertin and Lecu\'{e} (2008), who implement a two step procedure to reduce the dimensionality of a local estimate. The first step fits a local linear estimate with an $L_1$ or lasso type penalty, which identifies the relevant variables. This is followed by a second local linear fit using this reduced dimensionality. The lasso penalty they use is precisely the same as the third approach suggested in Section~\ref{se:approaches}.

We comment on the similarities and differences of these two approaches compared to the current presentation, which are summarised in Table~\ref{tabcomp}. Firstly the theoretical framework of the two other papers focus exclusively on the performance at a single point, while the LABAVS approach ensures uniformly oracle performance on the whole domain. The framework for the other two also assumes that variables are either active on the whole domain or redundant everywhere, while we have already discussed the usefulness of an approach that can adapt to variables that are redundant on various parts of the data. We believe this is particularly important, since local tests of variable significance will give the same results everywhere. Related to this, our method does not require an assumption of nonzero gradients (whether with respect to the bandwidth or variables) to obtain adequate theoretical performance, in contrast to the other methods. On the other hand, ensuring both uniform performance while allowing $d$ to be increasing is a quite challenging, so our presentation assumes $d$ is fixed,  while the others do not. It is also worth noting that the greedy approach of Lafferty and Wasserman potentially gives it an advantage in higher dimensional situations.

While all approaches work in a similar framework, the above discussion demonstrates that there are significant differences. Our methodology may be viewed as a generalisation of the work of Bertin and Lecu\'{e}, save for imposing fixed dimensionality. It can also be viewed as a competitor to the rodeo, and some numerical examples comparing the two are provided.

\begin{table}[htbp]\begin{center}
\caption{ Summary of locally adaptive bandwidth approaches  \label{tabcomp}}
\begin{tabular}{|l|c|c|c|} \hline
  &  LABAVS \T& Rodeo &  Bertin and \\
  &         &       &  Lecu\'{e} (2008) \B\\ \hline
Oracle performance on entire domain \T \B& \tick & \cross & \cross \\ \hline
Allows for locally redundant variables \T \B& \tick & \cross & \cross \\ \hline
Relevant variables allowed to have zero gradient \T \B& \tick & \cross & \cross \\ \hline
Theory allows dimension $d$ to increase with $n$ \T \B& \cross & \tick & \tick \\ \hline
Greedy algorithm applicable for higher dimensions \T \B& \cross & \tick & \cross \\ \hline
\end{tabular} \end{center} \end{table}

With regards to computation time, for estimation at a single point the rodeo is substantially faster, since calculating variable significance on a large grid of points is not required. If however we need to make predictions at a reasonable number of points, then Labavs is likely to be more efficient, since the grid calculations need only be done once, while rodeo requires a new set of bandwidth calculations for each point.

\section{Theoretical properties} \label{theory}

As mentioned in the introduction, a useful means of establishing the power of a model that includes variable selection is to compare it with an oracle model, where the redundant variables are removed before the modelling is undertaken. In the linear (and the parametric) context, we interpret the oracle property as satisfying two conditions as $n \to \infty$:
\begin{enumerate}
    \item the probability that the correct variables are selected converges to 1, and
    \item the nonzero parameters are estimated at the same asymptotic rate as they would be if the correct variables were known in advance.
\end{enumerate}
We wish to extend this notion of an oracle property to the nonparametric setting, where some predictors may be redundant. Here there are no parameters to estimate, so attention should instead be given to the error associated with estimating $g$. Below we define weak and strong forms of these oracle properties:
\begin{defn}
The {\it weak} oracle property in nonparametric regression is:
\begin{enumerate}
    \item the probability that the correct variables are selected converges to 1, and
    \item at each point $x$ the error of the estimator $\hg(x)$ decreases at the same asymptotic rate as it would if the correct variables were known in advance.
\end{enumerate}
\end{defn}
\begin{defn}
The {\it strong} oracle property in nonparametric regression is:
\begin{enumerate}
    \item the probability that the correct variables are selected converges to 1, and
    \item at each point $x$ the error of the estimator $\hg(x)$ has the same first-order asymptotic properties as it would if the correct variables were known in advance.
\end{enumerate}
\end{defn}
Observe that the weak oracle property achieves the correct rate of estimation while the strong version achieves both the correct rate and the same asymptotic distribution. The first definition is most analogous to its parametric counterpart, while the second is more ambitious in scope.

Here we establish the strong version of the nonparametric oracle property for the LABAVS algorithm, with technical details found in the appendix (Miller and Hall, 2010). We shall restrict attention to the case of fixed dimension. To apply to a situation of increasing dimension, we could add an asymptotically consistent screening method to reduce it back to fixed $d$. The treatment here focuses on local linear polynomials, partly for convenience but also recognising that the linear factors dominate higher order terms in the asymptotic local fit. Thus our initial fit is found by minimising the expression (\ref{eq:loclin}). We impose further conditions on the kernel $K$:
\beqn \begin{array}{p{12.7cm}}
$\int K(z)dz = 1$, $\int z^{(j)}K(z)dz = 0$ for each $j$,  $\int z^{(j)}z^{(k)} K(z)dz = 0$ when $j \neq k$ and $\int (z^{(j)})^2 K(z)dz = \mu_2(K)>0$, with  $\mu_2(K)$ independent of $j$. \label{eq:c4}
\end{array} \neqn
The useful quantity $R(K)$, depending on the choice of kernel, is defined as
$$ R(K) = \int K(z)^2 dz = \left\{ \int K^*(z^{(j)})^2 dz^{(j)}\right\}^d \,, $$
where $K^*$ is the univariate kernel introduced on page \pageref{kernel}. Let $a_n \asymp b_n$ denote the property that $a_n = O(b_n)$ and $b_n = O(a_n)$. We also require the following conditions (\ref{eq:cond}), needed to ensure uniform consistency of our estimators.
\beqn  \begin{array}{p{12cm}} \label{eq:cond}
\begin{enumerate}
\item  The support $\cC = \{x:f(x)>0\}$ of the random variable $X$ is compact. Further, $f$ and its first order partial derivatives are bounded and uniformly continuous on $\cC$, and $\inf_{x\in \cC} f(x) > 0$.
\item The kernel function $K$ is bounded with compact support and satisfies $|p(u) K(u) - p(v) K(v)| \le C_1||u-v||$ for some $C_1 > 0$ and all points $u,v$ in $\cC$. Here $p(u)$ denotes a single polynomial term of the form $\prod (u^{(j)})^{a_j}$ with the nonnegative integers $a_j$ satisfying $\sum a_j \le 4$. The bound $C_1$ should hold for all such choices of $p$.
\item The function $g$ has bounded and uniformly continuous partial derivatives up to order 2. If $(D^k g)(x)$ denotes the partial derivative
$$ \frac{\partial^{|k|}g(x)}{\partial (x^{(1)})^{k_1} \cdots \partial (x^{(d)})^{k_d}} \,,$$
with $|k| = \sum k_j$, then we assume that these derivatives satisfy, for some constant $C_2$,
$$ |h(u) - h(v)| \le C_2 || u - v || \,.$$
\end{enumerate}
\end{array} \neqn
$$ \begin{array}{p{12cm}}
\begin{enumerate}
\setcounter{enumi}{3}
\item $E(|Y|^\si) < \infty$ for some $\si > 2$.
\item The conditional density $f_{X|Y}(x|y)$ of $X_i$, conditional on $Y$, exists and is bounded.
\item $$ \mbox{For some $0<\xi<1$, } \frac{n^{1-2/\si} h^d}{\log n\{ \log n (\log\log n)^{1+\de}\}^{2/\si}} \to \infty \,.$$
\item The Hessian of $g$, $\cH_g$, is nonzero on a set of nonzero measure in $\cC$.
\end{enumerate}
\end{array} $$
The conditions in (\ref{eq:cond}), except perhaps the first, are fairly natural and not overly constrictive. For example, the sixth will occur naturally for any reasonable choice of $h$, while the second follows easily if $K$ has a bounded derivative. The last condition is purely for convenience in the asymptotics; if $\cH_g$ was zero almost everywhere then $g$ would be linear and there would be no bias in the estimate, improving accuracy. The first condition will not apply if the densities trail off to zero, rather than experiencing a sharp cutoff at the boundaries of $\cC$. However, in such circumstances our results apply to a subset of the entire domain, chosen so that the density did not fall below a specified minimum. Performance inside this region would then conform to the optimal accuracies presented, while estimation outside this region would be poorer. This distinction is unavoidable, since estimation in the tails is usually problematic and it would be unusual to guarantee uniformly good performance there.

Step 1 of the LABAVS Algorithm allows the initial bandwidth to be chosen globally or locally. Here we shall focus on the global case, where an initial bandwidth $H = \diag (h^2,\ldots,h^2)$ is used. Further, we assume that this $h$ is chosen to minimise the mean integrated squared error (MISE):
$$ E \left[ \int \{\hg(x) - g(x) \}^2 f(x) dx \right] \,,$$
where the outer expectation runs over the estimator $\hg$. It is possible to show that under our assumptions that
\beqn  h = \left[ \frac{d\si^2 R(K) A_\cC}{n \mu_2(K)^2 A_{\cH_g}} \right] ^{-1/(d+4)} \,. \label{eq:swan} \neqn
Notice in particular that $h \asymp n^{1/(d+4)}$. Details are given in Lemma~\ref{le:ace} in the appendix (Miller and Hall, 2010).

A key result in establishing good performance, in Theorem~\ref{th:1} below, is  uniform consistency of the local polynomial parameter estimates. It is a simplified version of a result by Masry (1996), and no proof is included.

\begin{thm} Suppose the conditions in (\ref{eq:cond}) hold and we use parameter estimates from a degree $p$ polynomial regression to estimate the partial derivatives of $g$. Then for each $k$ with $0 \le |k| \le p$ we have
$$ \sup_{x\in \cC} | (\widehat{D^k g})(x) - (D^k g)(x)| = O \left[ \left( \frac{\log n}{n h^{d+2|k|}} \right)^{1/2} \right] + O(h^{p-|k|+1}) \mbox{ almost surely.}$$ \label{th:1} \end{thm}

Since the partial derivative estimate at $x$ is proportional to the corresponding local polynomial coefficient, Theorem~1 ensures that the local polynomial coefficients are consistently estimated uniformly for suitable $h$. The scaling applied in (\ref{eq:scale}) does not impact on this, as the proof of Theorem~\ref{th:appr1} demonstrates.

Let $\cC^-$ denote the points $x\in\cC$ satisfying $\pa g(x)/\pa x^{(j)} = 0$ and $\pa^2 g(x)/\pa x^{(j)2} \ne 0$. That is, $\cC^-$ denotes the points where the true set of relevant variables changes. Notice that in the illustrative example in Section~\ref{se:eg} we had $\cC^- = \{ x \,|\, x^{(1)}=0 \,, x^{(2)}=0 \}$. The smoothness assumed of $g$ implies that $\cC^-$ has Lebesgue measure 0. Let $\de>0$ and let $\cO_\de$ be an open set in $\cC^-$ such that
\beqn \inf_{x \in \cC^+ \setminus \cO_\de,\, j \in \cA^+(x)} \ga_j = \de \,. \label{eq:duck} \neqn
Intuitively this means that on the set $\cC \setminus \cO_\de$ the relevant variables have the absolute value of their corresponding parameters $|\ga_j|$ bounded below by $\de>0$, while irrelevant variables have $\ga_j=0$. Thus we have a ``gap'' between the true and irrelevant variables in this region that we may exploit. The volume of $\cO_\de$ may be made arbitrarily small by choosing $\de$ small. Call the set $\chA^+(x)$ in the algorithm correct if the variables in it are the same as the set of variables $j$ with $\partial g(x)/\partial x^{(j)} \neq 0$. Denote the latter correct set by $\cA^+(x)$.

\begin{thm} Suppose $\de$ is given, $h$ is chosen to minimise squared error as in (\ref{eq:swan}), $\chA^+(x)$ is formed using the first approach in Section~\ref{se:approaches}, and $\la$ has a growth rate between arbitrary constant multiples of $h^2 (n\log n)^{1/2}$ and $hn^{1/2}$. If $f$ has bounded and uniformly continuous derivatives of degree 2, then the probability that $\chA^+(x)$ is correct on $\cC \setminus \cO_\de$ tends to 1 as $n \to \infty$. Furthermore, variables that are genuinely redundant everywhere will be correctly classified as such with probability tending to 1.
\label{th:appr1}  \end{thm}

The property (\ref{eq:duck}) ensures that the coefficients in the local linear fit are consistently estimated with error of order $O\{h(\log n)^{1/2}\}$. The adjustment in (\ref{eq:scale}) means that the actual coefficients estimated are of order $hn^{1/2}$ times this, so the range of $\la$ given is correct for separating true and redundant variables. The definition of $\cO_\de$ ensures that the classification is correct on $\cC \setminus \cO_\de$, while variables that are redundant everywhere will be recognised as such.

The next result ensures consistency for the second approach in Section~\ref{se:approaches}. We make one further assumption, concerning the error $\ep_i$.  Observe that this holds trivially if $\ep_i$ is bounded. Assume that:
\beqn \begin{array}{p{14cm}}
there exists $C_3$ such that $|E(\ep_i^\al)| \le C_3^{\al-2}\si^2$ for $\al=3,4,\ldots$.
\end{array} \label{eq:extra} \neqn

\begin{thm} Suppose $\de$ is given, $h$ is chosen to minimise squared error as in (\ref{eq:swan}), and $\chA^+(x)$ is formed using the second approach in Section~\ref{se:approaches}. Provided that $\la = o(h^2)$ and $h^4 \log n = o(\la)$, the probability that $\chA^+(x)$ is correct on $\cC \setminus \cO_\de$ tends to 1 as $n \to \infty$. Furthermore, variables that are genuinely redundant everywhere will be correctly classified as such with probability tending to 1.
\label{th:appr2} \end{thm}

The previous two results ensure that we have consistent variable selection for the first two approaches in Section~\ref{se:approaches}. Finally we can state and prove the strong oracle property for $\cC \setminus \cO_\de$. Although the result does not cover the whole space $\cC$, recall that we may make the area $\cO_\de$ arbitrarily small by decreasing $\de$. Furthermore, the proof implies that if we restricted attention to removing only those variables that are redundant everywhere, we would actually have the oracle property on the whole of $\cC$; however we sacrifice this performance on $\cO_\de$ to improve the fit elsewhere by adjusting for locally redundant variables. In the following theorem the matrix $\ddot{H}$ is the diagonal bandwidth matrix with bandwidth $\infty$ for globally redundant variables and $\ddot{h}$ for the other variables, where
$$ \ddot{h} = h \left( M(\ddot{H},H) \ddot{d}/d \right)^{1/4} \,.$$
Here $\ddot{d}$ denotes the number of variables that are not globally redundant.

\begin{thm} The estimates produced by the algorithm, where variable selection is performed using the first or second approach in Section~\ref{se:approaches}, satisfy the strong definition of the nonparametric oracle property on $\cC$. Further,  when there are locally redundant variables, squared estimation error is actually less that the oracle performance by a factor of $M[\{H^L(X),H^U(X)\},\ddot{H}] <1$. \label{th:big} \end{thm}

\section{Numerical properties} \label{se:num}

The examples presented in this section compare the performance of two versions of the LABAVS algorithm with ordinary least squares, a traditional local linear fit,  generalised additive models, tree-based gradient boosting and MARS. Table~\ref{tab:types} describes the approaches used. The implementations of the latter four methods were from the R packages {\sl locfit}, {\sl gam}, {\sl gbm} and {\sl polspline} respectively. Tuning parameters such as bandwidths for local methods, $\la$ in LABAVS, number of trees in boosting, and MARS model complexity, were chosen to give best performance for each method. The LABAVS models used the first variable selection approach of Section~\ref{se:approaches}. All the local methods used nearest neighbour bandwidths. The OLS linear model was included as a standard benchmark, but obviously will fail to adequately detect nonlinear features of a dataset.

\begin{table}[htbp]\begin{center}
\caption{ Approaches included in computational comparisons  \label{tab:types}}
\begin{tabular}{llll}
\hline
 Name &  Description \\
\hline
LABAVS-A & LABAVS with linear fit, all vars in final fit \\
LABAVS-B & LABAVS with linear fit, relevant vars only in final fit \\
LOC1      & Local linear regression     \\
OLS       & Ordinary least squares linear regression \\
GBM       & Boosting with trees, depth equal to three \\
GAM       & Generalised additive models with splines  \\
MARS      & Multivariate adaptive regression splines  \\
\hline
\end{tabular} \end{center} \end{table}

\begin{figure}[htb] \begin{center}
\includegraphics[width = 3.4in, height=3.4in]{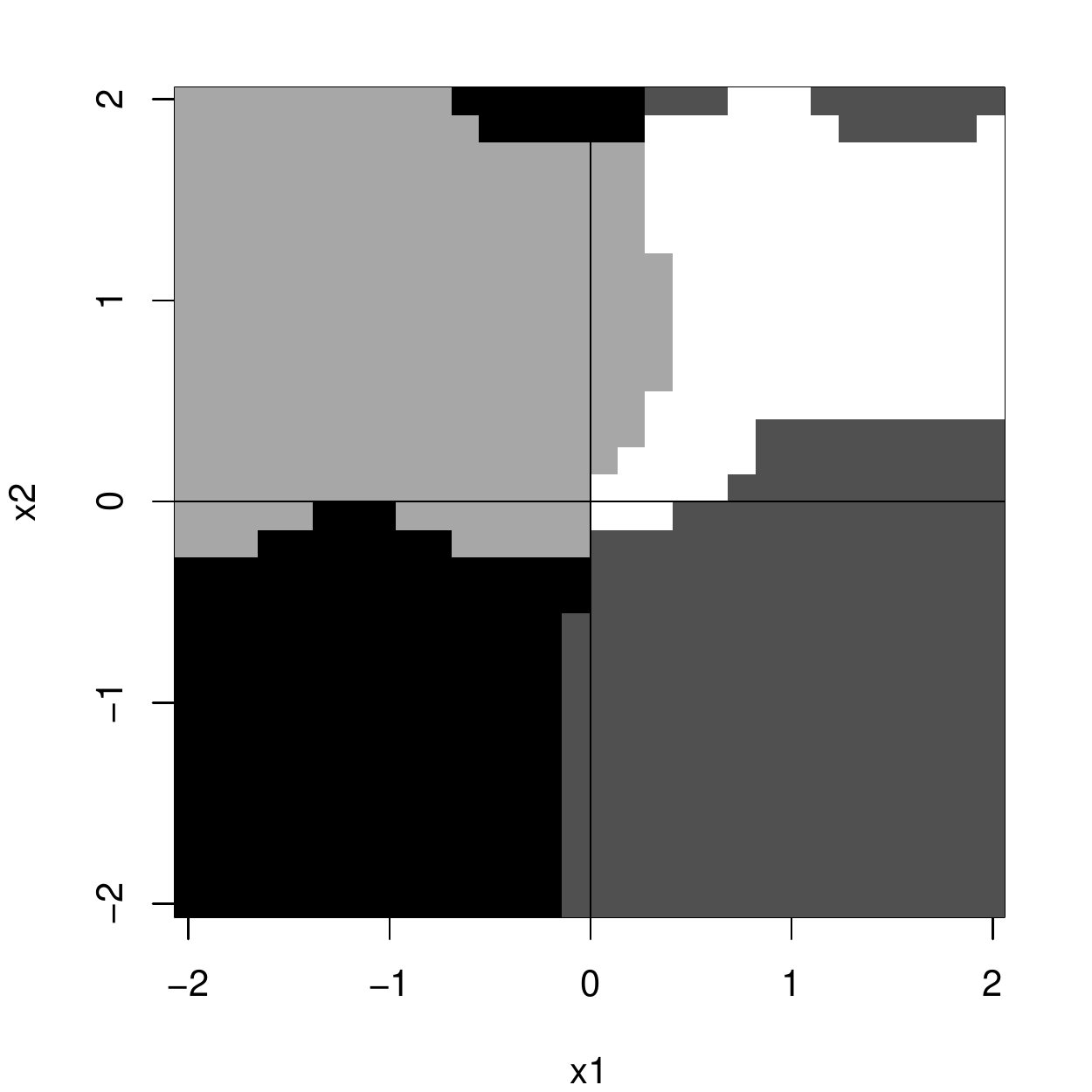}
\vspace{-8pt}
{\it \caption{ Plot of detected variable significance across subspace in Example 1. \label{fig3}}}
 \end{center} \end{figure}

\ni {\bf Example 1:} The example introduced in Section~\ref{se:eg} was simulated with $n=500$. The error for $Y_i$ was normal with standard deviation $0.3$. We first compare LABAVS to the rodeo and the methodology of Bertin and Lecu\'{e} (2008) at the four representative points in Figure~\ref{fig2}. Table~\ref{tabcomp1} shows the mean squared error of the prediction compared to the true value over 100 simulations. In all cases parameters were chosen to minimise this average error. At all points the LABAVS approach performed strongest. The method of Bertin and Lecu\'{e} (2008) performed poorly in situations where at least one variable is redundant; this is to be expected, since it excludes the variable completely and so will incorporate regions where it is actually important, causing significant bias. The rodeo also did not perform as well; we found it tended to overestimate the optimal bandwidths in redundant directions.

\begin{table}[htb]\begin{center}
\caption{ Mean squared prediction error on sample points in Example 1 \label{tabcomp1}}
\begin{tabular}{lcccc}
\hline
 Test Point & LABAVS-A & LABAVS-B & rodeo & Bertin and\\
&&&& Lecu\'{e} \\ \hline
(1,1)  & 0.0022 & 0.0022 & 0.0065 & 0.0023 \\
(1,-1) & 0.0011 & 0.0013 & 0.0015 & 0.0018 \\
(-1,1) & 0.0009 & 0.0011 & 0.0015 & 0.0013 \\
(-1,-1)& 0.0006 & 0.0007 & 0.0008 & 0.0013 \\
\hline
\end{tabular} \end{center} \end{table}

We then compared LABAVS with the other model approaches which are designed to make multiple predictions, rather than a specific point. For each simulation all the models were fitted and the average squared error was estimated using a separate test set of 500 observations. The simulation was run 100 times and the average error and its associated standard deviation for each model are recorded in Table~\ref{tab1}.

\begin{table}[htbp]\begin{center}
\caption{ Mean squared error sum of test dataset in Example 1 \label{tab1}}
\begin{tabular}{lcc}
\hline
\quad Approach &  Error & Std Dev\\
\hline
LABAVS-A  &    2.18   &       (0.71) \\
LABAVS-B  &    1.87   &       (0.65) \\
LOC1   &       2.31   &       (0.73) \\
OLS    &       42.85  &       (2.64) \\
GBM    &       2.47   &       (0.67) \\
GAM    &       5.93   &       (0.57) \\
MARS   &       2.35   &       (0.90) \\
\hline
\end{tabular} \end{center} \end{table}

Inspection of the results shows that the LABAVS models performed best, able to allow for the different dependencies on the variables. In particular the algorithm improved on the performance of the local linear model on which it is based. The local linear regression, the boosted model and MARS also performed reasonably, while the GAM struggled with the nonadditive nature of the problem, and a strict linear model is clearly unsuitable here.

To show how effectively variable selection is for  LABAVS, Figure~\ref{fig3} graphically represents the sets $\chA^+$ at each grid point for one of the simulations, with the darkest representing $\{\}$, the next darkest $\{1\}$, the next darkest $\{2\}$ and finally the lightest $\{1,2\}$. Here the variable selection has performed well; there is some encroachment of irrelevant variables into the wrong quadrants but the selection pattern is broadly correct. The encroachment is more prevalent near the boundaries since the bandwidths are slightly larger there, to cover the same number of neighbouring points.

\begin{figure}[tbhp] \begin{center}
\includegraphics[width = 5in, height=1.7in]{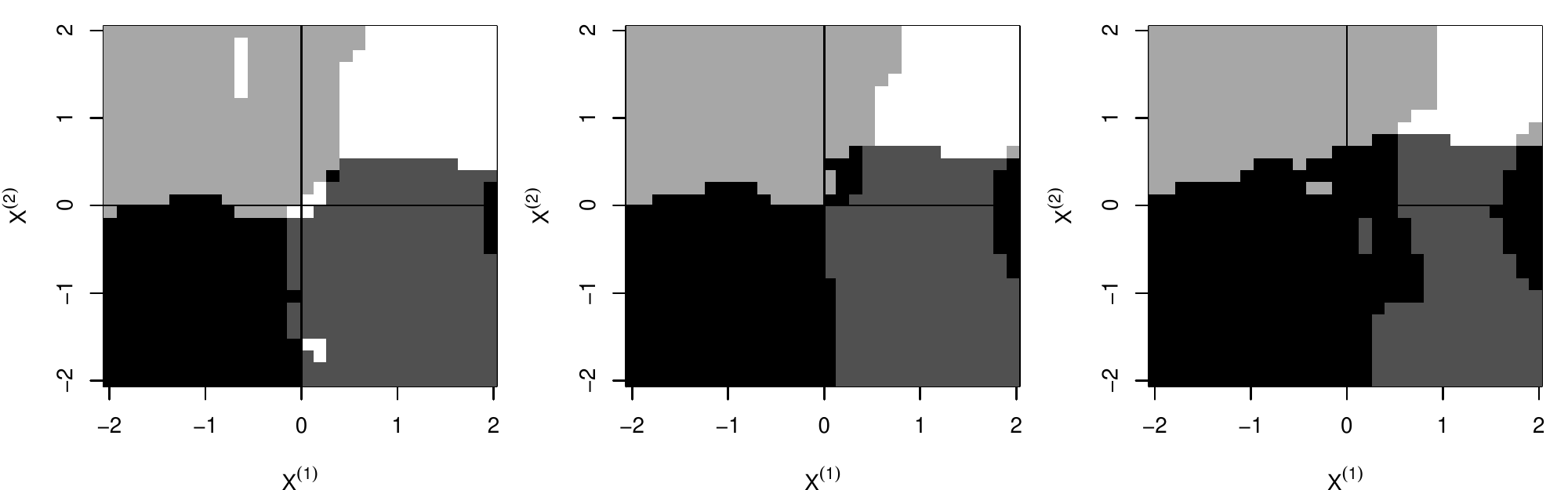}
\vspace{-8pt}
{\it \caption{ Plot of detected variable significance across subspace in Example 2, under various choices for $\la$. \label{figex2}}}
\end{center} \end{figure}

\ni {\bf Example 2:} We next show that LABAVS can effectively remove redundant variables completely. Retain the setup of Example~1, except that we add $d^*=d-2$ variables similarly distributed (uniform on [-2,2]), which have no influence on the response. Also, keep the parameters relating to the LABAVS fit the same as the previous example, except that the cutoff for hard threshold variable selection, $\la$ is permitted to vary. Table~\ref{tabex2} shows the proportion of times from 500 simulations that LABAVS effected complete removal of the redundant dimensions, for various $\la$ and $p^*$. Note that the cutoff level of 0.55 is that used in the previous example, and the two genuine variables were never completely removed in any of the simulations. The results suggest that to properly exclude redundant variables, a higher threshold is needed than would otherwise be the case. This causes the final model to be slightly underfit when compared to the oracle model, but this effect is not too severe; Figure~\ref{figex2} shows how the variable significance plots change for a particular simulation with different values of the cutoff. It is clear that the patterns are still broadly correct, and the results still represent a significant improvement to traditional linear regression.

\begin{table}[htbp]\begin{center}
\caption{ Proportion of simulations where redundant variables completely removed by LABAVS \label{tabex2}}
\begin{tabular}{lcccc}
\hline
\quad  & \multicolumn{4}{c}{Number of redundant dimensions}  \\
\hline
$\la$ & 1&2&3&4 \\ \hline
0.55& 0.394&0.086&0.034&0.038\\
0.65& 0.800&0.542&0.456&0.506\\
0.75& 0.952&0.892&0.874&0.864\\
0.85& 0.996&0.984&0.994&0.974\\
0.95& 0.998&1.000&1.000&0.992\\
\hline
\end{tabular} \end{center} \end{table}

\begin{figure}[tbh] \begin{center}
\includegraphics[width =5in, height=2.25in]{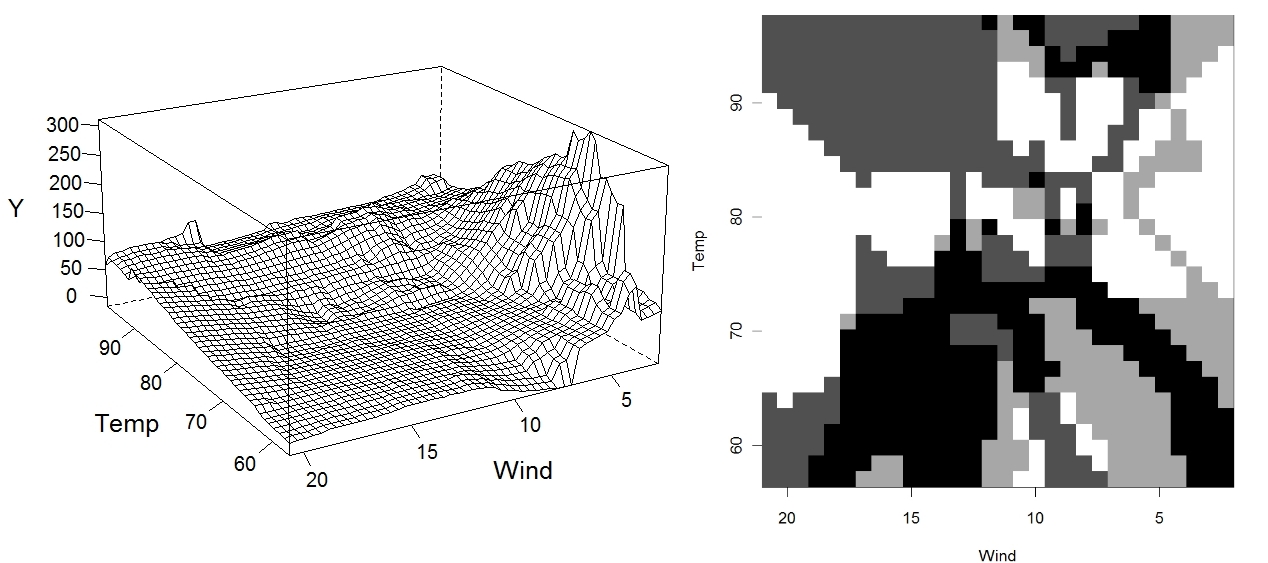}
\vspace{-8pt}
{\it \caption{ Ozone dataset smoothed perspective plot and variable selection plot. \label{fig4}}}
 \end{center} \end{figure}

\ni {\bf Example 3:} The first real data example used is the ozone dataset from Hastie~\etal (2001), p.175. It is the same as the air dataset in S-PLUS, up to a cube root transformation in the response. The dataset contains meteorological measurements for New York collected from May to September 1973. There are 111 observations in the dataset, a fairly moderate size. Our aim here is to predict the ozone concentration using two of the other variables, temperature and wind and scaled to unit variance when fitting the models. The smoothed perspective plot of the data in Figure~\ref{fig4} shows strong dependence on each of the two variables in some parts of the domain, but some sections appear flat in one or both directions in other parts. For example, the area surrounding a temperature of 70 and wind speed of 15 appears to be flat, implying that for reasonably low winds and high temperatures the ozone concentration is fairly stable. This suggests that LABAVS, by expanding the bandwidths here, could be potentially useful in reducing error. We performed a similar comparative analysis to that in Example 1, except that error rates were calculated using leave-one-out cross validation, where an estimate for each individual observations was made after using all other observations to build the model. The resulting mean squared errors and corresponding standard deviations are presented in Table~\ref{tab:Oz}. The results suggest that the data is best modelled using local linear methods, and that LABAVS offers a noticeable improvement over a traditional local fit, due to its ability to improve the estimate in the presence of redundant variables. The perspective plot in left panel of Figure~\ref{fig4} suggests a highly non-additive model, which may explain why GAM performs poorly. There is also a large amount of local curvature, which hinder the OLS, GBM and MARS fits. The right panel of Figure~\ref{fig4} shows the variable selection results for the linear version of LABAVS across the data support, using the same shading as in Figure~\ref{fig2}. We see that variable dependence is fairly complex, with all combinations of variables being significant in different regions. In particular, notice that the procedure has labelled both variables redundant in the region around $(70,15)$, confirming our initial suspicions. This plot is also highly suggestive, revealing further interesting features. For instance, there is also little dependence on wind when temperatures are relatively high. Such observations are noteworthy and potentially useful.

\begin{table}[htb]\begin{center}
\caption{ Cross-validated mean squared error sum for the ozone dataset  \label{tab:Oz}}
\begin{tabular}{lcc}
\hline
\quad Approach &  Error & Std Dev\\
\hline
LABAVS-A  &   277  &       (53) \\
LABAVS-B  &   284  &       (55) \\
LOC1      &   290  &       (55) \\
OLS       &   491  &       (110) \\
GBM       &   403  &       (118) \\
GAM       &   391  &       (98) \\
MARS      &   457  &       (115) \\
\hline
\end{tabular} \end{center} \end{table}

\ni {\bf Example 4:} As a second low-dimensional real data example, we use the ethanol dataset which has been studied extensively, for example by Loader (1999). The response is the amount of a certain set of pollutants emitted by an engine, with two predictors: the compression ratio of the engine and the equivalence ratio of air to petrol. There are 88 observations, a fairly moderate size. Inspection of the data shows strong dependence on the equivalence ratio, but the case for the compression ratio is less clear. This suggests LABAVS could be potentially useful in reducing error. We performed a similar analysis to that in Example 1, with the results are presented in Table~\ref{tab:Eth}.

\begin{table}[htb]\begin{center}
\caption{ Cross-validated mean squared error sum for the ethanol dataset  \label{tab:Eth}}
\begin{tabular}{lcc}
\hline
\quad Approach &  Error & Std Dev\\
\hline
LABAVS-A  &   0.075  &       (0.011) \\
LABAVS-B  &   0.085  &       (0.014) \\
LOC1      &   0.090  &       (0.012) \\
OLS       &   1.348  &       (0.128) \\
GBM       &   0.104  &       (0.020) \\
GAM       &   0.098  &       (0.012) \\
MARS      &   0.045  &       (0.008) \\
\hline
\end{tabular} \end{center} \end{table}

The results in Table~\ref{tab:Eth} show that this problem is particularly suited to MARS, which performed the best. After MARS, LABAVS produced the next strongest result, again improving on the traditional local linear model. The GBM and GAM models were inferior to the local linear fit.

\newpage
\appendix
\section{Technical details} \label{se:proofs}

We first prove the following lemma concerning the asymptotic behaviour of $h$.
\begin{lmm} The choice of $h$ that minimises the mean integrated squared error is asymptotically the minimiser of
\beqn (1/4)h^4\mu_2(K)^2 A_{\cH_g} + \si^2(nh^d)^{-1} R(K) A_\cC \,, \label{eq:sparrow} \neqn
where $R(K) = \int K(x)^2 dx$ for the function $K$, $A_{\cH_g} = \int \tr\{\cH_g(x)\}^2 f(x) dx$  and $A_\cC = \int_\cC 1 dx$. Further,
\beqn h = \left[ \frac{d\si^2 R(K) A_\cC}{n \mu_2(K)^2 A_{\cH_g}} \right] ^{1/(d+4)} \,. \label{eq:wolf} \neqn
\label{le:ace} \end{lmm}
\ni {\bf Proof:} Ruppert and Wand (1994) show that for $x$ in the interior of $\cC$ we have bias and variance expressions
$$ E\{\hg(x) - g(x)\} = (1/2) \mu_2(K) h^2\tr\{\cH_g(x)\} + o_P (h^2) \mbox{ , and} $$
$$ Var\{\hg(x)\} = \{ n^{-1}h^{-d} R(K) f(x)^{-1} \si^2\}\{1+o_P(1)\} \,. $$
Substituting these into the mean integrated squared error expression yields
\begin{eqnarray*}
MISE &=& \int E\{\hg(x) - g(x) \}^2 f(x) dx \\
    &=& \int \left[ E\{\hg(x)\}-g(x)]^2 + \Var\{\hg(x)\} \right] f(x) dx \\
    &=& \int (1/4) \mu_2(K)^2 h^4\tr\{\cH_g(x)\}^2 f(x) dx + o_P (h^4) \\
    &+& \int n^{-1}h^{-d} R(K) f(x)^{-1} \si^2 f(x) dx+ o_P(n^{-1}h^{-d}) \\
    &=& (1/4)h^4\mu_2(K)^2 A_{\cH_g} + \si^2(nh^d)^{-1} R(K) A_\cC + o_P(h^4 + n^{-1}h^{-d}) \,.
\end{eqnarray*}
This establishes the first part of the Lemma. Notice that assumptions (\ref{eq:c4}) and (\ref{eq:cond}) ensure that the factors $\mu_2(K)^2 A_{\cH_g} $ and $R(K) A_\cC$ are well defined and strictly positive.  Elementary calculus minimising (\ref{eq:sparrow}) with respect to $h$ completes the Lemma. $\blacksquare$ \\

Observe that we may express $Y_i$ using a first order Taylor expansion for $g$:
$$g(x) + \cD_g(x)^T(X_i-x) + \ep_i + T(x) \,,$$
where  the remainder term is $T(x) = \sum_{j,k} e_{j,k}(x)(X^{(j)}_i-x^{(j)})(X^{(k)}_i-x^{(k)})$ with the terms $e_{j,k}$ are uniformly bounded. For local linear regression we aim to show that our local linear approximation $\hga_0 + \hga^T (X_i-x)$ is a good approximation for this expansion and that the remainder behaves. The following two results are needed before proving the Theorem~\ref{th:appr1} and Theorem~\ref{th:appr2}. Firstly, the following version of Bernstein's Inequality may be found in Ibragimov and Linnik (1971), p169.
\begin{thm}[Bernstein's Inequality] Suppose $U_i$ are independent random variables, let $A_2=\sum_{i=1}^n \Var(U_i)$ and $S_n = \sum_{i=1}^n U_i$. Suppose further that for some $L>0$ we have
$$|E[\{U_i-E(U_i)\}^k]| \le \tfrac{1}{2} \Var(U_i) L^{k-2} k! \,.$$
Then
$$ P \{|S_n - E(S_n)| \ge 2t \sqrt{A_2} \} < 2e^{-t^2} \,. $$
\end{thm}

Secondly, the following lemma contains a proof which is applicable to many uniform convergence type results. The structure is similar to that of Masry (1996), although it is simplified considerably when using independent observations and Bernstein's Inequality. In the proof, let $C_4 = \sup f(x) < \infty$ and $C_5 = \inf f(x) >0$ for $x \in \cC$.

\begin{lmm} $$ \sup_{x \in \cC} | n^{-1} \tsum_i \ep_i K_H(X_i-x) | = O\left\{ \left( n^{-1}h^{-d} \log n   \right)^{1/2} \right\} \,.$$
\label{lm:Bern} \end{lmm}
\ni {\bf Proof:} Since $\ep_i$ is independent of $X_i$ and $E(\ep_i)=0$, we have $E\{\ep_i K_H(X_i-x) \} = 0$. As $\cC$ is compact we may cover it with $L(n) = (n / h^{d+2} \log n)^{d/2}$ cubes $I_1, \ldots, I_{L(n)}$, each with the same side length, proportional to $L(n)^{-1/d}$. Then
\begin{eqnarray*}
\sup_{x \in \cC} | n^{-1} \tsum_i \ep_i K_H(X_i-x) | &\le& \\
&& \hspace{-40pt} \max_m \sup_{x \in \cC \cap I_m} |n^{-1} \tsum\ep_i K_H(X_i-x) - n^{-1} \tsum \ep_i K_H(X_i-x_m) | \\
    && \hspace{-40pt} + \max_m |n^{-1} \tsum \ep_i K_H(X_i-x_m)| \,=\, Q_1 + Q_2
\end{eqnarray*}
From the second condition of (\ref{eq:cond}) we know that
\begin{eqnarray*}
|\ep_iK_H(X_i-x) - \ep_i K_H(X_i-x_m)| &\le& \frac{C_1\ep_i}{h^d} \para h^{-1}(x-x_m)  \para \\
&\le& \frac{C_1' \ep_i}{h^{d+1}} \left( \frac{h^{d+2} \log n}{n} \right)^{1/2} = C_1' \ep_i \left( \frac{\log n}{n h^d} \right)^{1/2}
\end{eqnarray*}
This expression is independent of $x$ and $m$, and so $Q_1 \le C_1' \left(\frac{\log n}{n h^d}\right)^{1/2} \left| n^{-1}  \sum \ep_i \right|$, which  implies that $Q_1 = O[(\frac{\log n}{n h^d} )^{1/2}]$. Now with regard to $Q_2$, notice that
\beqn P(Q_2 > \eta) \le L(n) \sup_x P\{|n^{-1} \tsum \ep_i K_H(X_i-x)|>\eta\}  \,. \label{eq:q2} \neqn
Letting $B_2 = \sup_x K(u)$ and using the first property in (\ref{eq:extra}) we see that for $k=3,4,\ldots$,
 \begin{eqnarray*}|E[\{\ep_i K_H(X_i-x)\}^\al ]| &\le& \si^2 C_3^{\al-2} \int K_H(u-x)^{\al}f(u)du \\
        &\le& \Var\{\ep_i K_H(X_i-x)\}(B_2C_3)^{\al-2} \,. \end{eqnarray*}
Also, if $B_3=\int K(u)^2 du$ we can show that $\Var\{K_H(X_i-x)\} \le C_4\si^2B_3h^{-d}$. We may let $n$ be large enough so that
$$ (B_4 \log n )^{1/2} \le \frac{\sqrt{\sum E\{\ep_i^2 K_H(X_i-x)^2\}}}{2B_2C_3} \,,$$
for some $B_4$ to be determined below. Then by Bernstein's inequality
\begin{eqnarray*}
P\bigg\{ && \hspace{-20pt} \Big| n^{-1} \sum \ep_i K_H(X_i-x) \Big| \ge 2(B_4\log n)^{1/2} \Big(\frac{\si^2C_4 B_3^2}{nh^d}\Big)^{1/2} \bigg\} \\
&\le& P\left\{ | \sum \ep_i K_H(X_i-x)| \ge 2(B_4\log n)^{1/2} \sqrt{\sum E(\ep_i^2 K_H(X_i-x)^2)} \right\} \\
&\le& 2e^{-B_4\log n} \le \frac{2}{n^{-B_4}} \,.
\end{eqnarray*}
Comparing this inequality to (\ref{eq:q2}) and choosing $B_4$ large enough then the expression $2L(n)n^{-B_4}$ is summable, by the Borel-Cantelli lemma we may conclude that $Q_2 = O[(\frac{\log n}{n h^d} )^{1/2}]$ and the lemma is proved. $\blacksquare$ \\

In a similar fashion it is also possible to prove, letting $Z_i = X_i - x$ and $\ze = (n^{-1}h^{-d}\log n)^{1/2}$,
\begin{eqnarray}
&\sup_x | n^{-1}\sum_i K_H(Z_i)- E\{K_H(Z_i)\}| &= O(\ze) \label{eq:magpie} \\
&\sup_x | n^{-1}\sum_i K_H(Z_i)^2- E\{K_H(Z_i)^2\}| &= O(h^{-d} \ze)  \\
&\sup_x | n^{-1}\sum_i Z_i^{(j)}K_H(Z_i)- E\{Z_{i}^{(j)}K_H(Z_i)\}| &= O(h \ze) \\
&\sup_x | n^{-1}\sum_i \ep_i Z_i^{(j)}K_H(Z_i)- E\{\ep_i Z_{i}^{(j)}K_H(Z_i)\}| &= O(h \ze) \\
&\sup_x | n^{-1}\sum_i Z_i^{(j)}Z_i^{(k)}K_H(Z_i)- E\{Z_{i}^{(j)}Z_{i}^{(k)}K_H(Z_i)\}| &= O(h^2 \ze) \label{eq:moose} \\
&\sup_x | n^{-1}\sum_i e_{jk} Z_i^{(j)}Z_i^{(k)}K_H(Z_i)- E\{e_{jk} Z_{i}^{(j)}Z_{i}^{(k)}K_H(Z_i)\}| &= O(h^2 \ze) \\
&\sup_x | n^{-1}\sum_i e_{jk} Z_i^{(j)}Z_i^{(k)}Z_i^{(l)}K_H(Z_i) \qquad\qquad& \label{eq:robin}\\
&\qquad\qquad\qquad\qquad- E\{e_{jk} Z_{i}^{(j)}Z_{i}^{(k)}Z_i^{(l)}K_H(Z_i)\}| &= O(h^3 \ze) \nonumber
\end{eqnarray}
Standard treatment of the expectation integrals reveals that
\begin{eqnarray}
E\{K_H(Z_i)\} &=& f(x)+ O(h) \label{eq:raven}\\
E\{K_H(Z_i)^2\} &=& h^{-d}\{f(x)R(K)+ O(h)\} \\
E\{Z_i^{(j)} K_H(Z_i)\} &=& O(h^2) \\
E\{\ep_i Z_i^{(j)} K_H(Z_i)\} &=& 0 \label{eq:epfun} \\
E\{Z_i^{(j)}Z_i^{(k)} K_H(Z_i)\} &=& O(h^2) \\
E\{e_{jk} Z_i^{(j)}Z_i^{(k)} K_H(Z_i)\} &=& O(h^2) \\
E\{e_{jk} Z_i^{(j)}Z_i^{(k)}Z_i^{(l)} K_H(Z_i)\} &=& O(h^4) \label{eq:crow}
\end{eqnarray}

If $h\asymp n^{-1/(d+4)}$, as it will be under Lemma~\ref{le:ace}, then the asymptotic rates in the expectations (\ref{eq:raven})-(\ref{eq:crow}) will dominate those of the deviations (\ref{eq:magpie})-(\ref{eq:robin}), with the execption of (\ref{eq:epfun}). We may then conclude that, uniformly on $x$,
\begin{eqnarray}
n^{-1}\sum_i K_H(Z_i) &=& f(x) + O(h) \label{eq:C1} \\
n^{-1}\sum_i K_H(Z_i)^2 &=& h^{-d}\{f(x)R(K) + O(h)\} \label{eq:Cb} \\
n^{-1}\sum_i \ep_i K_H(Z_i) &=& O(h)\\
n^{-1}\sum_i Z_i^{(j)}K_H(Z_i) &=& O(h^2) \label{eq:Cc} \\
n^{-1}\sum_i \ep_i Z_i^{(j)}K_H(Z_i) &=& O(h^2) \\
n^{-1}\sum_i Z_i^{(j)}Z_i^{(k)}K_H(Z_i) &=& O(h^2) \\
n^{-1}\sum_i e_{jk}Z_i^{(j)}Z_i^{(k)}K_H(Z_i) &=& O(h^2) \\
n^{-1}\sum_i e_{jk}Z_i^{(j)}Z_i^{(k)}Z_i^{(l)}K_H(Z_i) &=& O(h^4) \label{eq:C2}
\end{eqnarray}

\ \\
{\ni \bf Proof of Theorem~\ref{th:appr1}:} From Lemma~\ref{le:ace} we know that an estimator of $h$ that minimises mean integrated squared error will satisfy $h\asymp n^{-1/(d+4)}$. Theorem~\ref{th:1} then implies that
$$ \sup_{x \in \cC, j = 1,\ldots,d} | (\widehat{D^j g})(x) - (D^j g)(x) | = O(h \sqrt{\log n}) \,.$$
Notice that the estimates $\hga_j$ at $x$ in the minimisation (\ref{eq:loclin}) are exactly the estimates $(\widehat{D^j g})(x)$. The adjusted parameter estimates $\hbe_j$ in (\ref{eq:appr1}) therefore satisfy
\beqn \hbe_j = (\widehat{D^j g})(x) \{ \tsum (X_i^{(j)} - \bX_x^{(j)})^2 K_H(X_i-x)\}^{1/2} \,. \label{eq:hick} \neqn
Let $\be_j = (D^j g)(x)\{nh^2 \mu_2(K) f(x) \}^{1/2}$. We aim to show that $\hbe$ converges to $\be$ sufficiently fast uniformly in $x$.
\begin{eqnarray}
\sup_{x\in \cC,\, j} |\hbe_j - \be_j| &\le& \sup_{x,j} \big| \{nh^2 \mu_2(K) f(x)\}^{1/2} \{(\widehat{D^j g})(x) - (D^j g)(x)\} \big| \nonumber \\
 && \hspace{-60pt} + \sup_{x,j} \big| (\widehat{D^j g})(x) [ \{ \tsum (X_i^{(j)}-\bX_x^{(j)})^2 K_H(X_i-x)\}^{1/2} - \{ nh^2 \mu_2(K) f(x) \}^{1/2} ] \big| \nonumber \\
&\le& O(h^2\sqrt{n\log n}) \nonumber \\
 && \hspace{-60pt} + A_1 \sup_{x,j} \big| \{ \tsum (X_i^{(j)}-\bX_x^{(j)})^2 K_H(X_i-x)\}^{1/2} - \{ nh^2 \mu_2(K) f(x) \}^{1/2} \big|
  \label{eq:cow}
\end{eqnarray}
In the first term of the last line we use the fact that $(D^j g)$ is bounded and $(\widehat{D^j g})$ converges uniformly so may be bound be some constant $A_1$, and for the second term we use the boundedness of $f(x)$ and (\ref{eq:hick}).

Focusing on the first term, note that
$$\sup_{x,\, j} |\bX_x^{(j)} - x^{(j)} | = \sup_{x,\, j} \Big| \frac{\sum(X_i^{(j)}-x^{(j)})K_H(X_i-x)}{\sum K_H(X_i-x)} \Big| = O(h^2) \,,$$
using (\ref{eq:C1}) and (\ref{eq:Cc}). Thus
\begin{eqnarray*}
\sum_i (X_i-\bX_x^{(j)})^2 K_H(X_i-x) &=& \sum \{ X_i^{(j)} - x^{(j)} + O(h^2) \}^2 K_H(X_i-x) \\
    &=& O(nh^4) + \sum (X_i^{(j)} - x^{(j)})^2 K_H(X_i-x) \,,
\end{eqnarray*}
again using (\ref{eq:C1}) and (\ref{eq:Cc}). Now we consider the expectation of $(X_i^{(j)} - x^{(j)})^2 K_H(X_i-x)$ carefully,
\begin{eqnarray*}
E\{(X_i^{(j)} - x^{(j)})^2 K_H(X_i-x)\} &=& \int (u^{(j)} - x^{(j)})^2K_H(u-x)f(u)du \\
&=& h^2\int (z^{(j)})^2 K(z) f(x+hz) dz \\
&=& h^2\int (z^{(j)})^2 K(z) \{f(x)+h z^T \cD_f(x) + O(h^2)\} dz \\
&=& h^2 \mu_2(K) f(x) + O(h^4) \,.
\end{eqnarray*}
The differentiability assumptions in the statement of the Theorem ensure that this formulation is uniform over all $x$ in $\cC$. Using this and (\ref{eq:moose}) in (\ref{eq:cow}) and noting that if $x \to 0$ then $(1+x)^{1/2}-1 = O(x)$, we see that
\begin{eqnarray*}
\sup_{x\in \cC,\, j} |\hbe_j - \be_j| &\le& A_1 \sup_{x,j} | \{ nh^2 \mu_2(K) f(x) + O(n h^2 \sqrt{\log n}))\}^{1/2} - \{ nh^2 \mu_2(K) f(x) \}^{1/2} | \\
&+& O(h^2\sqrt{n\log n})  \\
&=& \sup_x A_1 \{ nh^2 \mu_2(K)f(x)\}^{1/2} | \{ 1+ O(h^2 \sqrt{\log n}) \}^{1/2} - 1 | \\
&+& O(h^2\sqrt{n\log n}) \\
&=& \sup_x A_1 \{ nh^2 \mu_2(K)f(x)\}^{1/2} O(h^2 \sqrt{\log n}) + O(h^2\sqrt{n\log n}) \\
&=& O( h^2 \sqrt{n\log n})
\end{eqnarray*}
We do not need to worry about small nonzero values of $\be_j$ by our assumption on $\cO_\de$, so the nonzero $\be_j$ grow at $O(n^{1/2}h)$. Further, the estimate error of $\hbe_j$ is $O( h^2 \sqrt{n\log n})$ uniformly in $x$. A $\la$ that grows at some rate between these two as, suggested in the Theorem, will be able to separate the true variables from the redundant ones with probability tending to 1. $\blacksquare$

\ \\
{\ni \bf Proof of Theorem~\ref{th:appr2}:} Let $\hmu_0$, $\hmu$ be the parameter estimates for the case where the $j$th variable is removed from consideration, so $\mu^{(j)}=0$. Theorem~\ref{th:1} ensures that the maximum distance for the estimators $\hga_0, \hmu_0$ from $g(x)$ is $O(\ze)$ and similarly $\hga, \hmu$ converge to the derivative $\cD_g(x)$ at $O(\ze h^{-1})$, with the exception of $\hmu^{(j)}=0$. Thus we may expand the sum of squares difference and use results (\ref{eq:C1})--(\ref{eq:C2}):
\begin{eqnarray*}
SS_j(x) \hspace{-10pt} &-\hspace{-10pt}& SS(x) \\
&=& n^{-1}\sum\{Y_i - \hmu_0 - \hmu^T Z_i \}^2K_H(Z_i) - n^{-1}\sum\{Y_i - \hga_0 - \hga^T Z_i \}^2K_H(Z_i) \\
&=& n^{-1}\sum K_H(Z_i)\bigg[ \{O(\ze) + \ep_i + \cD_g(x)^{(j)}Z_i^{(j)}+T(x) + O(\ze h^{-1})\tsum_k Z_i^{(k)} \}^2 \\
&-&  \{O(\ze) +\ep_i +T(x) + O(\ze h^{-1})\tsum_k Z_i^{(k)} \}^2 \bigg] \\
&=& n^{-1} \sum K_H(Z_i) \bigg[ O(\ze^2) + \ep_i O(\ze) + T(x) O(\ze) + O(\ze^2 h^{-1})\tsum_k Z_i^{(k)} \\
&+& O(\ze) \cD_g(x)^{(j)}Z_i^{(j)} + 2\ep_i \cD_g(x)^{(j)}Z_i^{(j)} + 2T(x) \cD_g(x)^{(j)}Z_i^{(j)} \\
&+& O(\ze h^{-1}) \cD_g(x)^{(j)}\tsum_k Z_i^{(j)} Z_i^{(k)} + (\cD_g(x)^{(j)}Z_i^{(j)})^2 + O(\ze h^{-1}) \ep_i \tsum_k Z_i^{(k)} \\
&+& O(\ze h^{-1}) T(x) \tsum_k Z_i^{(k)} + O(\ze^2 h^{-2})\tsum_{k,\, l} Z_i^{(k)} Z_i^{(l)}\bigg] \\
&=& O(\ze^2) + (\cD_g^{(j)})^2 O(h^2) \,.
\end{eqnarray*}
This shows the behaviour of the numerator in our expression. Note that our assumption on $\cO_\de$ ensures that when $|\cD_g^{(j)}|$ is nonzero it is bounded away from 0, so true separation is possible. In a similar fashion to that above we may expand and deal with the denominator $n^{-1}\sum(Y_i - \hga_0 - \hga^T Z_i )^2 K_H(Z_i)$. The dominating term here is the asymptotic expectation  of $n^{-1}\sum\ep_i^2 K_H(Z_i)$, which tends to $\si^2 f(x)$, and everything else converges to zero at $h$ or faster, uniformly in $x$. Therefore, so as long as $\la$ shrinks faster than $h^2$ but slower than $\ze^2 = h^4 \log n$, the variable selection will be uniformly consistent. $\blacksquare$ \\

Before proving Theorem~\ref{th:big}, we prove the following three lemmas. The first allows us to separate out the effects of various variables in the LABAVS procedure. The latter two are concerned with the change in estimation error for local and global variable redundance respectively.

\begin{lmm} Let $\cB_1$ and $\cB_2$ be disjoint subsets of $\{1,\ldots,d\}$ such that $\cB_1 \cup \cB_2 = \{1,\ldots,d\}$. The final estimates of the LABAVS procedure would be the same as applying the bandwidth adjustment, that is steps 3 and 4, in the procedure twice; the first time only expanding the bandwidths at $x$ of those variables in $\chA^-(x) \cap \cB_1$ to the edges of the maximal rectangle and shrinking those remaining, and the second time expanding the variables in $\chA^-(x) \cap \cB_2$ and shrinking the variables in $\chA^+(x)$.
\label{lm:big1} \end{lmm}
\ni {\bf Proof:} Choose $x \in \cC$. With some slight abuse of notation, since the bandwidths are possibly asymmetric, let $H_1(x)$ denote the adjusted bandwidths after the first step of the two-step procedure, with shrunken variables having bandwidth $h_1$. Similarly let $H_2(x)$ denote the bandwidths after the second step, with bandwidth on the shrunken variables $h_2$. Further, let $d_1(x)$ equal the cardinality of $\chA^+(x) \cup \cB_1$ and $d_2(x)$ equal the cardinality of $\chA^+(x)$. The bandwidths for the redundant variables are expanded to the edges of the maximal rectangle, so we need only show that the resulting shrunken bandwidth is the same as when applying the one-step version of the algorithm. Using expression (\ref{newfac2}) we know that
\begin{eqnarray*}
h_1  &=& M(H_1, H) = h \left[ \frac{E\{d_1(X)\} E\{ V(X,H_1)\}}{d E\{ V(X,H)\}} \right] \,, \mbox{ and} \\
h_2  &=& M(H_2, H_1) = h_1 \left[ \frac{E\{d_2(X)\} E\{ V(X,H_2)\}}{E\{d_1(X)\} E\{ V(X,H_1)\}} \right] \,.
\end{eqnarray*}
Substituting the first expression into the second gives
$$ h_2 = h \left[ \frac{E\{d_2(X)\} E\{ V(X,H_2)\}}{d E\{ V(X,H)\}} \right] \,, $$
which recovers the equation in the one-step bandwidth adjustment. Thus the bandwidths are unchanged for every $x \in \cC$. $\blacksquare$

\begin{lmm} Suppose that $h$ is chosen to minimise squared error as in (\ref{eq:swan}). Also, suppose that the LABAVS procedure identifies that no variables are globally redundant but some (possibly all) variables are locally redundant and that the local redundancy takes place on a set of non-zero measure.  Then the LABAVS procedure reduces the overall MISE of the estimation of $g$ by a factor of $M[\{H^L(X),H^U(X)\},H] < 1$.
\label{lm:big2} \end{lmm}
\ni {\bf Proof:} We shall ignore the difficulties associated with incorrect selection on $\cO_\de$, as it only affects an arbitrarily small subset of the domain. With probability tending to one we have correct variable classification, so we work under this assumption. Since some variables are relevant in some regions, the choice of $h'$ is well defined. Pick $x \in \cC$ and let $u^+$ denote the components of $d$-vector $u$ indexed by $\cA^+(x)$ and $u^-$ be the residual components. We can express the density $f(x)$ as $f(x^+,x^-)$ so the relevant and redundant components may be treated separately. From (\ref{eq:C1}) and (\ref{eq:Cb}) we know that
$$ V(x,H) = \frac{h^{-d} \{ R(K)f(x) + O(h)\}}{\{ f(x) + O(h)\}^2} = h^{-d} \left\{ \frac{R(k)}{ f(x)}+O(h)\right\} \,.$$
Taking an expectation over $x$ we see that
\beqn E\{V(X,H)\} = h^{-d} R(K) A_\cC + O(h^{-(d-1)}) \,. \label{eq:mongoose} \neqn
For convenience let $H^*(x)$ denote the asymmetric bandwidths $H^L(x)$ and $H^U(x)$.  We now show that the factor $M\{H^*(X),H\}$ is less than 1. Firstly observe that $V(x,H) = V(x, H^*(x))$ whenever $\cA^+(x) = (1,\ldots,d)$. Consider the case when $\cA^+(x) \neq (1,\ldots,d)$. In particular assume that $k$ components are redundant at $x$. We see that
\begin{eqnarray*}
E\{ K_{H^*} (X_i-x)^2 \} &=& \int \int \bigg\{(h')^{-(d-k)} f(x^{+}+h'z^{+}, u^{-}) \prod_{j\in \cA^+(x)}^n K^*(z^{(j)})^2  \\
&& \cdot  \prod_{j\in \cA^-(x)}^n h_j^*(x)^{-2}K^*\{h^*_j(x)^{-1}(u^{(j)}-x^{(j)})\}^2 \bigg\}dz^{+}du^{-} \\
&=& (h')^{-(d-k)} R(K)^{(d-k)/d} \bigg[ O(h')  \\
&& \hspace{-30pt} +\int \prod_{j\in \cA^-(x)}^n h_j^*(x)^{-2} K^*\{h^*_j(x)^{-1}(u^{(j)}-x^{(j)})\}^2 f(x^{+}, u^{-})du^{-}  \bigg]\\
&=& (h')^{-(d-k)}\{B_1(x) + O(h') \} \,,
\end{eqnarray*}
where $B_1(x)$ is a uniformly bounded and strictly positive number depending only on $x$. An argument using Bernstein's Theorem similar to that in Lemma~\ref{lm:Bern} shows that the uniform bound of $n^{-1} \sum K_{H^*}(X_i-x)^2$ away from $E\{ K_{H^*} (X_i-x)^2 \}$ is $O[(h')^{-(d-k)} \{n(h')^{(d-k)}\}^{-1/2}\sqrt{\log n} ]$, so we may deduce that
$$ n^{-1} \sum K_{H^*}(X_i-x)^2 = (h')^{-(d-k)} \{ B_1(x) + O(h') \} \,.$$
In a similar fashion we can show that
$$ n^{-1} \sum K_{H^*}(X_i-x) = f_{X^{+}}(x^{+})B_2(x) + O(h') \,,$$
where $B_2(x)$ is a uniformly bounded and strictly positive number depending only on $x$. This leads to
\beqn  V(x,H^*) = (h')^{-(d-k)}\left\{ \frac{B_1(x)}{B_2(x)^2} + O(h') \right\} \,. \label{eq:dove} \neqn
Let $\cE$ denote the event that $\cA^+(X)=(1,\ldots,n)$ and $\cE^c$ the complement. We know $P(\cE^c)>0$ by assumption and also that given $\cE^c$ is true, $V(X,H^*) = O\{(h')^{-(d-1)}\}$ from (\ref{eq:dove}). Thus as $n\to \infty$, for some strictly positive constants $B_3$, $B_4$ and $B_5$,
\begin{eqnarray*}
\frac{E\{ V(X,H^*)\} }{E\{ V(X,H)\}} &=& \frac{ P(\cE ) E\{ V(X,H^*)\ | \cE\} + P(\cE^c ) E\{ V(X,H^*)\ | \cE^c\} }{ P(\cE) E\{ V(X,H)\ | \cE\} + P(\cE^c ) E\{ V(X,H)\ | \cE^c\} } \\
&=& \frac{h^{-d}\{B_3+O(h)\} + O\{(h')^{-(d-1)}\}}{h^{-d}\{B_3+O(h)\} + h^{-d}\{B_4+O(h)\}} \,,
\end{eqnarray*}
where $B_3$, $B_4$ are constants satisfying $B_3 \ge 0$ and $B_4 > 0$. But from our definition of $h'$ in (\ref{newfac2}) and the definition of $M(H^*(X),H)$, we may deduce that
$$ \left(\frac{h'}{h}\right)^4 \asymp \frac{B_3 + O\{h^d(h')^{-(d-1)}\}}{B_3 + B_4} + O(h) \,.$$
From this expression it follows that both sides must be less than 1 in the limit. Thus we have $M(H^*,H) <1$ asymptotically, as required.

Ruppert and Wand (1994) show that for a point $x$ in the interior of $\cC$, using a bandwidth matrix $H$, $\Var(\hat{g}(x)|X_1, \ldots, X_n)$ is equal to
$$\si^2 e_1^T (\dX^TW\dX)^{-1} \dX^T W^2 \dX(\dX^TW\dX)^{-1}e_1 \,,$$
 where $\dX$ is the $n \times (p+1)$ matrix $({\bf 1},X)$, $e_1$ is a $p$-vector with first entry $1$ and the others $0$, and $W$ is an $n \times n$ diagonal matrix with entries $K_H(X_i-x)$. This variance may be reexpressed as
$$ \si^2 \frac{\tsum K_H(X_i-x)^2}{\{\tsum K_H(X_i-x)\}^2} e_1^T (\dX^T\dX)^{-1} \dX^T \dX(\dX^T\dX)^{-1}e_1\,.$$
Taking ratios of the expectations for the variance factors under the adjusted and initial bandwidths recovers the expression $M(H^*,H)$ in (\ref{eq:fac}). Thus the variance term in the MISE is reduced by a factor of $M(H^*,H)$. Furthermore, the bias term in the MISE, in which we may ignore the zero bias contributed be the $n$th variable where it is redundant, is reduced by a factor of $(h'/h)^4$ which, from (\ref{newfac2}), is strictly less than the factor $M(H^*,H)$. Thus the MISE is reduced by the factor $M(H^*,H)$ as required. $\blacksquare$

\begin{lmm} Suppose that $h$ is chosen to minimise squared error as in (\ref{eq:swan}). Also, suppose that the LABAVS procedure finds that all variables are relevant everywhere in $\cC$ except for a single variable $X^{(j)}$, which is globally irrelevant. Then the LABAVS procedure reduces the overall MISE of the estimation of $g$ by a factor of $M(H^*,H) < 1$. Furthermore the resulting bandwidth $h'$ is asymptotically optimal, in the sense that it minimises the $d-1$ dimensional MISE expression.
\label{lm:big3} \end{lmm}
\ni {\bf Proof:} Let $\cC'$ denote the $d-1$ dimensional space formed by removing the irrelevant variable and denote the volume of this space by $A_{\cC'}$. We know that our initial $h$ satisfies (\ref{eq:wolf}). By similar reasoning it follows that we are required to show that our adjusted bandwidth is asymptotically equal to
\beqn h_{\mbox{\small opt}} = \left[ \frac{(d-1)\si^2 R(K)^{(d-1)/d} A_{\cC'}}{n \mu_2(K)^2 A_{\cH_g}} \right] ^{1/(d+3)} \,, \label{eq:sheep} \neqn
which is the bandwidth the minimises MISE in the reduced dimension case. Here $A_{\cC'}$ denotes the volume of the $d-1$ dimensional case. Equivalently, combining (\ref{eq:wolf}) and (\ref{eq:sheep}), it is sufficient to show in the limit that
\beqn \frac{(h')^{d+3}}{h^{d+4}} = \frac{(d-1) A_\cC'}{d R(K)^{1/d} A_\cC} \,. \label{eq:goat} \neqn
Arguments similar to those in the previous Lemma can be made to show
$$ n^{-1} \sum K_{H^*}(X_i-x)^2 = (h')^{-(d-1)} \{ R(K)^{(d-1)/d} f_{X^{(-n)}}(x^{(-n)})+O(h') \} \mbox{, and}$$
$$ n^{-1} \sum K_{H^*}(X_i-x) =  \{ R(K)^{(d-1)/d} f_{X^{(-n)}}(x^{(-n)})+O(h') \,.$$
Thus
\begin{eqnarray*}
E\{ V(X,H^*) \} &=& (h')^{-(d-1)} \bigg\{  R(K)^{(d-1)/d} \int \int f_{X^{(-n)}}(x^{(-n)})^{-1} f_{X^{(-n)}}(x^{(-n)}) \\ && \cdot f_{X^{(n)}|X^{(-n)}}(u^{(n)}) du^{(n)} du^{(-n)} + O(h') \bigg\} \\
&=&  (h')^{-(d-1)} \{ A_{\cC'} R(K)^{(d-1)/d} + O(h') \}
\end{eqnarray*}
Combining this with (\ref{eq:mongoose}) and (\ref{newfac2}) gives
$$ \left(\frac{h'}{h} \right)^4 = \frac{d-1}{d} M(H^*,H) = \frac{d-1}{d} \frac{h^d}{(h')^{d-1}} \left\{ \frac{A_{\cC'}}{A_\cC} R(K)^{-1/d} + O(h) \right\} \,.$$
Rearranging this last expression and letting $n\to \infty$ leads to the required expression (\ref{eq:goat}). Note that (\ref{eq:goat}) also implies that $(h'/h)^{d+3} h^{-1}$ is asymptotically constant, so $h'/h \to 0$. This in turn implies that $(h'/h)^4 = M(H^*,H)(d-1)/d$ tends to zero so asymptotically $M(H^*,H)<1$ as required. The argument that the MISE is in fact reduced by the factor $M(H^*,H)$ is entirely analogous to the previous Lemma. $\blacksquare$

\ \\
{\ni \bf Proof of Theorem~\ref{th:big}:} Correct variable selection at every point $x\in \cC$ with probability tending to 1 on the set $\cC \setminus \cO_\de$ for locally redundant variables, and $\cC$ for globally redundant variables, is guaranteed by Theorem~\ref{th:appr1} or Theorem~\ref{th:appr2}. For a given point $x$, repeated application of Lemma~\ref{lm:big1} allows us to consider the eventually result by adjusting the bandwidths for any partition of variables in any order. Choose an order in which globally redundant variables are treated first, one at a time, followed by a final adjustment for those variables that are locally redundant. Lemma~\ref{lm:big3} ensures that when allowing for each globally redundant variable, the resulting bandwidths in the remaining variables is asymptotically optimal. This means that the strong nonparametric oracle property is satisfied after the global bandwidth adjustments. Lemma~\ref{lm:big2} provides the quantification of the additional benefit resulting from the local variable removal. $\blacksquare$

\end{document}